\documentclass[preprint,3p]{elsarticle}
\usepackage{amsmath}
\usepackage{amssymb}
\usepackage{caption} 
 \usepackage[english]{babel} 
 \usepackage{float} 
 \usepackage{graphics} 
 \usepackage{graphicx} 
\usepackage[utf8]{inputenc}
\usepackage{mathtools}
\usepackage{capt-of}
\usepackage{calrsfs}
\usepackage{lipsum}  
 \usepackage{epsfig}
\newtheorem{definition}{Definition}[section]
\newtheorem{prop}{Proposition}[section]
\newtheorem{pro}{Property}[section]
\newdefinition{remark}{Remark}[section]
\newdefinition{fe}{Feature}[section]
\newproof{proof}{Proof}
\newproof{pot}{Proof of Theorem \ref{thm2}}
\newproof{poot}{Proof of Corollary \ref{co1}}
\numberwithin{equation}{section}
\newdefinition{exe}{Example}[section]
\journal{}
\begin{document}
\begin{frontmatter}

\title{Fractional Sturm-Liouville eigenvalue problems, I }
\author{
M. Dehghan\footnote[1]{dehghan@math.carleton.ca} and A. B. Mingarelli\footnote[2]{angelo@math.carleton.ca}}

\address{$^{1}$Department of Mathematics, Sari Branch, Islamic Azad University, Sari, Iran.\\
$^{2}$School of Mathematics and Statistics, Carleton University, Ottawa, Canada}

\begin{abstract}
We introduce and present the general solution of three two-term fractional differential equations of mixed Caputo/Riemann Liouville  type. We then solve a Dirichlet type Sturm-Liouville eigenvalue problem for a fractional differential equation derived from a special  composition of a Caputo and a Riemann-Liouville operator on a finite interval where the boundary conditions are induced by evaluating Riemann-Liouville integrals at those end-points. For each $1/2<\alpha<1$ it is shown that there is a finite number of real  eigenvalues, an infinite number of non-real eigenvalues, that the number of such real eigenvalues grows without bound as $\alpha \to 1^-$, and that the fractional operator converges to an ordinary two term Sturm-Liouville operator as $\alpha \to 1^-$ with Dirichlet boundary conditions.  Finally, two-sided estimates as to their location are provided as is their asymptotic behavior as a function of $\alpha$.
\end{abstract}

\end{frontmatter}
 \section{Introduction}
Fractional Sturm-Liouville Problems (FSLP)  generalize classical SLP in that the ordinary derivatives are replaced by {\it fractional derivatives}, or derivatives of fractional order. As an introduction the interested reader may wish to consult the great variety of works on the subject, starting with books such as  \cite{kil, pod, sam}. It turns out that such fractional equations have seen many applications in science and engineering in problems such as capturing the dynamical behaviors of amorphous materials, e. g., polymer and porous media \cite{metzler,rossik} or the superior modeling of the anomalous diffusion in materials with memory, e. g., viscoelastic materials for which the mean square variance grows faster (superdiffusion) or slower (subdiffusion) than in a Gaussian process \cite{bouch,met,shl}. 
 
One of the major differences between fractional and ordinary derivatives lies in the global nature of the former and to the local nature of the latter.  Thus in order to get information on the fractional derivative of a function at a given point one needs to have a knowledge of the original function on a semi-interval! The ordinary circular trigonometric functions that occur in the classical theory are now to be replaced by Mittag-Leffler functions.  So, it is conceivable and practical that physical problems involving past memory and also equations with delay be governed using such  fractional derivatives \cite{poos,sing,maina}.
 
The question at the core of this paper involves the determination of analytical solutions and subsequent spectrum of a fractional Sturm-Liouville equation with zero potential function and fixed end homogeneous boundary conditions.  As can be seen below even such an innocuous question can require a tremendous effort. To the best of our knowledge, the vast majority of scientific papers have paid  attention to a numerical solution of such equations (see \cite{tomas, zayer, mad1, mad2}).
 
Of note here is Klimek and Agrawal's work \cite{kliag1} in which the fractional operators are different from those usually defined in the literature in the sense that they contain left and right Riemann-Liouville and left and right Caputo fractional derivatives. They were also able to prove that, under suitable conditions, if the eigenvalues exist then they are real and that the eigenfunctions enjoy the same orthogonal property as in the classical case \cite{kliag2}.

 \section{Preliminaries}  
 \subsection{Fractional Calculus}
We recall some basic definitions in ODE's and results of importance in the fractional calculus. 
The Wronskian of two differentiable functions $f(t)$  and $g(t)$ is defined as usual by
$\mathcal{W}(f,g)(t)=f(t)g'(t)-f'(t)g(t).$  If $f(t)$ and $g(t)$ are two solutions of a second order linear differential equation and $\mathcal{W}(f,g)(t)\neq 0$ for some $t$, then we say that $f(t)$ and $g(t)$ are linearly independent and these functions are referred to as a fundamental set of solutions (FSS).
 \begin{definition}
The left and the right Riemann-Liouville fractional integrals $\mathcal{I}^{\alpha}_{a^+}$ and $\mathcal{I}^{\alpha}_{b^-}$ of order $\alpha\in \mathbb{R}^+$ are defined by
\begin{equation}\label{lfi}
\mathcal{I}^{\alpha}_{a^+}f(t):=\frac{1}{\Gamma(\alpha)}\int_a^t\frac{f(s)}{(t-s)^{1-\alpha}}ds,\quad t \in (a,b],
\end{equation}
and
\begin{equation}\label{rfi}
\mathcal{I}^{\alpha}_{b^-}f(t):=\frac{1}{\Gamma(\alpha)}\int_t^b\frac{f(s)}{(s-t)^{1-\alpha}}ds,\quad t \in [a,b),
\end{equation}
respectively. Here $\Gamma (\alpha)$ denotes  Euler's Gamma function. The following property is easily verified.
 \end{definition}
 \begin{pro}\label{c}
\normalfont 
 For a constant $C$, we have $\mathcal{I}_{a^+}^{\alpha}C=\frac{(t-a)^{\alpha}}{\Gamma(\alpha+1)}\cdot C$.
 \end{pro}
  \begin{definition}
The left and the right Caputo fractional derivatives $^{c}\mathcal{D}^{\alpha}_{a^+}$ and $^{c}\mathcal{D}^{\alpha}_{b^-}$ are defined by
\begin{equation}\label{lcfd}
^{c}\mathcal{D}^{\alpha}_{a^+}f(t):=\mathcal{I}^{n-\alpha}_{a^+}\circ\mathcal{D}^nf(t)=\frac{1}{\Gamma(n-\alpha)}\int_a^{t}\frac{f^{(n)}(s)}{(t-s)^{\alpha-n+1}}ds, \quad t>a,
\end{equation}
and 
\begin{equation}\label{rcfd}
^{c}\mathcal{D}^{\alpha}_{b^-}f(t):=(-1)^n\mathcal{I}^{n-\alpha}_{b^-}\circ\mathcal{D}^nf(t)=\frac{(-1)^n}{\Gamma(n-\alpha)}\int_t^{b}\frac{f^{(n)}(s)}{(s-t)^{\alpha-n+1}}ds, \quad t<b,
\end{equation}
respectively, where $f$ is sufficiently differentiable and $n-1\leq \alpha<n$.
 \end{definition}
  \begin{definition}
Similarly, the left and the right Riemann-Liouville  fractional derivatives $\mathcal{D}^{\alpha}_{a^+}$ and $\mathcal{D}^{\alpha}_{b^-}$ are defined by
\begin{equation}\label{lrlfd}
\mathcal{D}^{\alpha}_{a^+}f(t):=\mathcal{D}^n\circ\mathcal{I}^{n-\alpha}_{a^+}f(t)=\frac{1}{\Gamma(n-\alpha)}\frac{d^n}{dt^n}\int_a^{t}\frac{f(s)}{(t-s)^{\alpha-n+1}}ds, \quad t>a,
\end{equation}
and 
\begin{equation}\label{rrlfd}
\mathcal{D}^{\alpha}_{b^-}f(t):=(-1)^n\mathcal{D}^n\circ\mathcal{I}^{n-\alpha}_{b^-}f(t)=\frac{(-1)^n}{\Gamma(n-\alpha)}\frac{d^n}{dt^n}\int_t^{b}\frac{f(s)}{(s-t)^{\alpha-n+1}}ds, \quad t<b,
\end{equation}
respectively, where $f$ is sufficiently differentiable and $n-1\leq \alpha<n$.
 \end{definition}
 \begin{pro}
\normalfont
 (Lemma 2.19 \cite{klim}.) Assume that $0<\alpha<1$, $f\in AC[a,b]$ and $g\in L^p(a,b)(1\leq p \leq \infty)$. Then the following fractional integration by parts formula holds
 \begin{equation}\label{part}
 \int_a^bf(t)\mathcal{D}^{\alpha}_{a^+}g(t)dt=\int_a^b g(t) ^{c}\mathcal{D}^{\alpha}_{b^-}f(t)dt+f(t)\mathcal{I}^{1-\alpha}_{a^+}g(t)|_{t=a}^{t=b}.
 \end{equation}
 \end{pro}
\begin{pro}
\normalfont
(Property 2.1 \cite{kil}.) Let $0<\alpha<\beta$, then the following identities hold:
\begin{equation}\nonumber
\begin{split}
\mathcal{I}^{\alpha}_{a^+}(t-a)^{\beta-1}&=\frac{\Gamma(\beta)}{\Gamma(\beta+\alpha)}(t-a)^{\beta+\alpha-1},\\
\mathcal{D}^{\alpha}_{a^+}(t-a)^{\beta-1}&=\frac{\Gamma(\beta)}{\Gamma(\beta-\alpha)}(t-a)^{\beta-\alpha-1},\\
\mathcal{I}^{\alpha}_{b^-}(b-t)^{\beta-1}&=\frac{\Gamma(\beta)}{\Gamma(\beta+\alpha)}(b-t)^{\beta+\alpha-1},\\
\mathcal{D}^{\alpha}_{b^-}(b-t)^{\beta-1}&=\frac{\Gamma(\beta)}{\Gamma(\beta-\alpha)}(b-t)^{\beta-\alpha-1}.
\end{split}
\end{equation}
\end{pro} 
 \begin{pro}\label{betaf} 
\normalfont
For $a \leq t < b$, $0 < \alpha < 1$ we have,
\begin{eqnarray*}
\mathcal{I}_{a^+}^{\alpha}\left ((b-t)^{\alpha-1}\right )&=& - \frac{(b-t)^{2\alpha -1}}{\Gamma(\alpha)\,}\, \int_{\frac{b-a}{b-t}}^1 (w-1)^{\alpha-1}\, w^{\alpha -1}\, dw,\\
&=& \frac{(b-t)^{2\alpha -1}}{\Gamma(\alpha)\,}\, \left (B \left ( \frac{b-a}{b-t};\alpha,\alpha \right) - B(1;\alpha,\alpha)\right ).
\end{eqnarray*}
where $B(z;\gamma,\theta)$ is the ``Incomplete Beta function" defined by 
\[
B(z;\gamma,\theta)=\int_0^z u^{\gamma-1}(1-u)^{\theta-1}du.
\]
 \end{pro}
 \begin{pro}
\normalfont
 (Lemma 2.4 \cite{kil}.) If $\alpha>0$ and $f\in L^p(a,b)(1\leq p\leq\infty)$, then the following equalities 
 \begin{equation}\nonumber
\begin{split}
\mathcal{D}^{\alpha}_{a^+}\circ\mathcal{I}^{\alpha}_{a^+}f(t)=f(t),\\
\mathcal{D}^{\alpha}_{b^-}\circ\mathcal{I}^{\alpha}_{b^-}f(t)=f(t),
\end{split}
\end{equation}
hold almost everywhere on $[a,b]$.
 \end{pro}
  \begin{pro}\label{p4}
\normalfont
 (Lemma 2.5 and Lemma 2.6 \cite{kil}.) If $0<\alpha<1$, $f\in L^1(a,b)$ and $\mathcal{I}^{1-\alpha}_{a^+}f,\mathcal{I}^{1-\alpha}_{b^-}f\in AC[a,b]$, then the following equalities 
 \begin{equation}\nonumber
\begin{split}
\mathcal{I}^{\alpha}_{a^+}\circ\mathcal{D}^{\alpha}_{a^+}f(t)=f(t)-\frac{(t-a)^{\alpha-1}}{\Gamma(\alpha)}\mathcal{I}^{1-\alpha}_{a^+}f(t)|_{t=a},\\
\mathcal{I}^{\alpha}_{b^-}\circ\mathcal{D}^{\alpha}_{b^-}f(t)=f(t)-\frac{(b-t)^{\alpha-1}}{\Gamma(\alpha)}\mathcal{I}^{1-\alpha}_{b^-}f(t)|_{t=b},
\end{split}
\end{equation}
hold a.e. on $[a,b]$.
 \end{pro}
 \begin{pro}
\normalfont
 (Lemma 2.21 \cite{kil}.) Let $\Re(\alpha)>0$ and $f(t)\in L^{\infty}(a,b)$ or $f(t)\in C[a,b]$. If $\Re(\alpha)\notin\mathbb{N}$ or $\alpha \in \mathbb{N}$, then 
 \begin{equation}\nonumber
\begin{split}
^{c}\mathcal{D}^{\alpha}_{a^+}\circ\mathcal{I}^{\alpha}_{a^+}f(t)=f(t),\\
^{c}\mathcal{D}^{\alpha}_{b^-}\circ\mathcal{I}^{\alpha}_{b^-}f(t)=f(t).
\end{split}
\end{equation}
\end{pro}
 
  \begin{pro}\label{p6}
\normalfont
 (Lemma 2.22 \cite{kil}.) Let $0<\alpha\leq 1$. If $f\in AC[a,b]$, then 
 \begin{equation}\nonumber
\begin{split}
\mathcal{I}^{\alpha}_{a^+}\circ     ^{c}\mathcal{D}^{\alpha}_{a^+} f(t)=f(t)-f(a),\\
\mathcal{I}^{\alpha}_{b^-}\circ      ^{c}\mathcal{D}^{\alpha}_{b^-}f(t)=f(t)-f(b).
\end{split}
\end{equation}
\end{pro}

 \subsection{The Mittag-Leffler function}
 The function $E_{\delta}(z)$ defined by 
 \begin{equation}\label{mittag1}
 E_{\delta}(z):=\sum_{k=0}^{\infty}\frac{z^{k}}{\Gamma(\delta k+1)},\quad (z\in \mathbb{C}, \Re(\delta)>0),
 \end{equation}
 was introduced by Mittag-Leffler \cite{mit}. In particular, when $\delta=1$ and $\delta=2$, we have
 \begin{equation}\label{e1}
 E_1(z)=e^z,\qquad E_2(z)=\cosh(\sqrt{z}).
\end{equation}
 The Mittag-Leffler function $E_{\delta,\theta}(z)$, generalizing the one in (\ref{mittag1}) is normally defined by
 \begin{equation}\label{mittag2}
 E_{\delta,\theta}(z):=\sum_{k=0}^{\infty}\frac{z^k}{\Gamma(\delta k+\theta)},\quad (z,\theta\in \mathbb{C}, \Re(\delta)>0).
 \end{equation}
 Of course, when $\theta=1$, $E_{\delta,\theta}(z)$ coincides with the Mittag-Leffler function (\ref{mittag1}):
 \begin{equation}\label{e2}
 E_{\delta,1}(z)=E_{\delta}(z).
 \end{equation}
 Two other particular cases of (\ref{mittag2}) are as follows:
 \begin{equation}\label{pmi}
 E_{1,2}(z)=\frac{e^z-1}{z},\quad E_{2,2}(z)=\frac{\sinh(\sqrt{z})}{\sqrt{z}}.
 \end{equation}
 \begin{pro}\label{asym}
\normalfont
If $0<\delta<2$ and $\mu\in(\frac{\delta \pi}{2},\min(\pi,\delta \pi))$, then function $E_{\delta,\theta}(z)$ has the following asymptotic expansion as $|z|\rightarrow\infty$ (see \cite{kil})
 \begin{eqnarray}
E_{\delta,\theta}(z)=\left\{
\begin{array}{ll}
\frac{1}{\delta}z^{\frac{1-\theta}{\delta}}\exp(z^{\frac{1}{\delta}})-\sum_{k=1}^N\frac{1}{\Gamma(\theta-\delta k)}\frac{1}{z^k}+O(\frac{1}{z^{N+1}}), \qquad |\arg(z)|\leq \mu, \\
-\sum_{k=1}^N\frac{1}{\Gamma(\theta-\delta k)}\frac{1}{z^k}+O(\frac{1}{z^{N+1}}), \qquad \mu\leq|\arg(z)|\leq \pi.
\end{array}\right.
 \label{asy}
\end{eqnarray}
\end{pro}
 \subsection{The Laplace transform}
 \begin{definition}
 The Laplace transform of a function $f(t)$ defined for all $t\geq 0$, is the function $F(s)$ defined by
 \[
F(s)=\mathfrak{L}\{f(t)\}:=\int_0^{\infty}e^{-st}f(t)dt,
 \]
 whenever the integral exists, where $s$ is the frequency parameter.
 \end{definition}
  \begin{definition}
 The inverse Laplace transform of a function $F(s)$ is then given by the line integral
 \[
 f(t)=\mathfrak{L}\{F(s)\}:=\frac{1}{2\pi i}\lim_{T\rightarrow\infty}\int_{\gamma-iT}^{\gamma+iT}e^{st}F(s)ds
 \] 
 where the integration is done along the vertical line $\Re(s)=\gamma$ in the complex plane such that $\gamma$ is greater than the real part of all the singularities of $F(s)$.
 \end{definition}

 \begin{definition}
The convolution of $f(t)$ and $g(t)$ supported on only $[0,\infty)$ is defined by
 \[
 (f\ast g)(t)=\int_0^tf(s)g(t-s)ds,\qquad f,g:[0,\infty)\rightarrow \mathbb{R}.
 \]
whenever the integral exists.
 \end{definition}
\begin{pro}
\normalfont
For $\Re(q)>-1$, then 
\begin{gather}\label{tq}
\mathfrak{L}\{t^q\}=\frac{\Gamma(q+1)}{s^{q+1}}\\
\mathfrak{L}^{-1}\{s^q\}=\frac{1}{t^{q+1}\Gamma(q)}\label{ils}
 \end{gather}
\end{pro}
\begin{pro}\label{lapder}
\normalfont If $f(t)$ is assumed to be a differentiable function and its derivative is of exponential type, then $\mathfrak{L}\{f'(t)\}=s\mathfrak{L}\{f(t)\}-f(0)$.
 \end{pro}
\begin{pro}\label{con}
\normalfont
$\mathfrak{L}\{(f\ast g )(t)\}=\mathfrak{L}\{f(t)\}\cdot  \mathfrak{L}\{g(t)\}$.
 \end{pro}
\begin{pro}\label{sal}
\normalfont
By definition of the left fractional integral (\ref{lfi}), we can rewrite
\[
\mathcal{I}_{0^+}^{\alpha}f(t)=\frac{1}{\Gamma(\alpha)}(f(t)\ast\frac{1}{t^{1-\alpha}}).
\]
So, by virtue of \eqref{tq} and Property~\ref{con}, we have
\begin{equation}\nonumber
\begin{split}
\mathfrak{L}\{\mathcal{I}_{0^+}^{\alpha}f(t)\}&=\frac{1}{\Gamma(\alpha)}\mathfrak{L}\{f(t)\}\cdot \mathfrak{L}\{\frac{1}{t^{1-\alpha}}\}\\
&=\frac{1}{s^{\alpha}}\mathfrak{L}\{f(t)\}.
\end{split}
\end{equation}
\end{pro}
\begin{pro}\label{lc}
The definition of the left Caputo fractional derivative (\ref{lcfd}) and Properties~\ref{sal} and \ref{lapder}, gives us, for $0<\alpha<1$,
\begin{equation}\nonumber
\begin{split}
\mathfrak{L}\{^c\mathcal{D}_{0^+}^{\alpha}f(t)\}&=\mathfrak{L}\{\mathcal{I}_{0^+}^{1-\alpha}\mathcal{D}f(t)\}\\
 &=\frac{1}{s^{1-\alpha}}\mathfrak{L}\{\mathcal{D}f(t)\}\\
 &=\frac{1}{s^{1-\alpha}}\left(s\mathfrak{L}\{f(t)\}-f(0)\right)\\
 &=s^{\alpha}\mathfrak{L}\{f(t)\}-s^{\alpha-1}f(0).
\end{split}
\end{equation}
\end{pro}
\begin{pro}\label{lrl}
Similarly, the definition of the left Riemann-Liouville fractional derivative (\ref{lrlfd}),  and Properties~\ref{sal} and \ref{lapder} gives us, for $0<\alpha<1$, 
\begin{equation}\nonumber
\begin{split}
\mathfrak{L}\{\mathcal{D}_{0^+}^{\alpha}f(t)\}&=\mathfrak{L}\{\mathcal{D}\mathcal{I}_{0^+}^{1-\alpha}f(t)\}\\
 &=s\mathfrak{L}\{\mathcal{I}_{0^+}^{1-\alpha}f(t)\}-\mathcal{I}_{0^+}^{1-\alpha}f(t)|_{t=0}\\
 &=s\left(\frac{1}{s^{1-\alpha}}\mathfrak{L}\{f(t)\}\right)-\mathcal{I}_{0^+}^{1-\alpha}f(t)|_{t=0}\\
 &=s^{\alpha}\mathfrak{L}\{f(t)\}-\mathcal{I}_{0^+}^{1-\alpha}f(t)|_{t=0}.
\end{split}
\end{equation}
\end{pro}
\section{Fundamental set of Solutions (FSS)}
In this section, we consider three features of a differently defined fractional Sturm-Liouville operator. This operator is a composition of a left Riemann-Liouville fractional derivative with a right Caputo fractional derivative as follows:
\begin{equation}\label{fe1}
^{c}\mathcal{D}_{b^-}^{\alpha}\circ\mathcal{D}^{\alpha}_{a^+}y(t)=0,\qquad 0<\alpha<1.
\end{equation}
By a solution of \eqref{fe1} is meant a function $y \in AC[a,b]$ such that ${D}^{\alpha}_{a^+} \in AC[a,b]$.
Applying the right fractional integral on (\ref{fe1}) and using the second relation of Property~\ref{p6} we obtain
\[
\mathcal{D}^{\alpha}_{a^+}y(t)-\mathcal{D}^{\alpha}_{a^+}y(t)|_{t=b}=0.
\]
Now, by taking the left fractional integral of the above equation followed by the first relation in Property~\ref{p4}, and then Property~\ref{c}, we get
\begin{equation}\label{sf1}
y(t)=\frac{(t-a)^{\alpha-1}}{\Gamma(\alpha)}\cdot \mathcal{I}_{a^+}^{1-\alpha}y(t)|_{t=a}+\frac{(t-a)^{\alpha}}{\Gamma(\alpha+1)}\cdot \mathcal{D}_{a^+}^{\alpha}y(t)|_{t=b},
\end{equation}
in which $\mathcal{I}_{a^+}^{1-\alpha}y(t)|_{t=a}$ and $\mathcal{D}_{a^+}^{\alpha}y(t)|_{t=b}$ are constants to be determined by imposing one or more initial/boundary conditions. The two solutions $y_1(t)=\frac{(t-a)^{\alpha-1}}{\Gamma(\alpha)}$ and $y_2(t)=\frac{(t-a)^{\alpha}}{\Gamma(\alpha+1)}$ are then fundamental solutions of  (\ref{fe1}) since both of them satisfy the equation (\ref{fe1}) separately, as is readily verified, they are linearly independent and finally
\[
\mathcal{W}(y_1,y_2)(t)=\frac{1}{\alpha}\frac{(t-a)^{2\alpha-2}}{\Gamma^2(\alpha)},
\]
 is not identically zero in the whole interval $[a,b]$ and has just discontinuity at point $a$. Observe that $y'_2(t)=y_1(t)$.
\begin{remark}
It is worth noting that as $\alpha$ approaches $1$, (\ref{fe1}) reduces to $-y^{\prime\prime}=0$, our two solutions $y_1(t)$ and $y_2(t)$ converge to $1$ and $(t-a)$ respectively, and (\ref{sf1}) then becomes $y(t)=y(a)+y'(b)(t-a)$, 
which authenticates classical results. We also have $\lim_{\alpha\rightarrow 1}\mathcal{W}(y_1,y_2)(t)=1$.
\end{remark}
Associated to \eqref{fe1} is another similar but quite different composition i.e.,
\begin{equation}\label{fe2}
\mathcal{D}_{b^-}^{\alpha}\circ   ^{c}\mathcal{D}^{\alpha}_{a^+}y(t)=0,\qquad 0<\alpha<1.
\end{equation}
Applying the right fractional integral on (\ref{fe1}) and using the second relation in Property~\ref{p4}, we have
\[
^{c}\mathcal{D}^{\alpha}_{a^+}y(t)-\frac{(b-t)^{\alpha-1}}{\Gamma(\alpha)}\cdot \mathcal{I}_{b^-}^{1-\alpha}\circ   ^{c}\mathcal{D}^{\alpha}_{a^+}y(t)|_{t=b}=0.
\]
Now, by taking the left fractional integral on the above equation,  using  Properties~\ref{p6} and \ref{betaf}, and introducing the function $\psi(t;a,b,\alpha) $ by
\[
\psi(t;a,b,\alpha)=\frac{(b-t)^{2\alpha-1}}{\Gamma^2(\alpha)}\cdot \left \{B \left (\frac{b-a}{b-t};\alpha,\alpha\right ) - B(1;\alpha,\alpha)\right \}
\]
we obtain
\begin{equation}\label{sf2}
y(t)=y(a)+\mathcal{I}_{b^-}^{1-\alpha}\circ   ^{c}\mathcal{D}^{\alpha}_{a^+}y(t)|_{t=b}\cdot\psi(t;a,b,\alpha).
\end{equation}
So, we can say FSS is $\{1,\psi(t;a,b,\alpha)\}$, since both of them satisfy the equation (\ref{fe2}) separately and their Wronskian is not identically zero in $[a,b]$,  having discontinuities at $a$ and $b$.
\begin{pro}\label{xik}
It is obvious that $\psi(a;a,b,\alpha)=0$ and if $\frac{1}{2}<\alpha<1$ we have, by a simple application of both Leibniz's  and L'Hospital's Rule,
\[
\lim_{t\rightarrow b}\psi(t;a,b,\alpha)=\frac{(b-a)^{2\alpha-1}}{(2\alpha-1)\Gamma^2(\alpha)}.
\]
The function $\psi(t;a,b,\alpha)$ has a discontinuity at $t=b$ when $0<\alpha\leq \frac{1}{2}$.
\end{pro}
In order to determine the constant value $\mathcal{I}_{b^-}^{1-\alpha}\circ   ^{c}\mathcal{D}^{\alpha}_{a^+}y(t)|_{t=b}$, we substitute the value $t=b$ into equation (\ref{sf2}) and by virtue of Property~\ref{xik} for $\frac{1}{2}<\alpha<1$, we get
\begin{equation}\label{ax}
\mathcal{I}_{b^-}^{1-\alpha}\circ   ^{c}\mathcal{D}^{\alpha}_{a^+}y(t)|_{t=b}=\frac{\left(y(b)-y(a)\right)(2\alpha-1)\Gamma^2(\alpha)}{(b-a)^{2\alpha-1}}.
\end{equation}
Next, we substitute the right hand side of (\ref{ax}) into (\ref{sf2}) and we obtain the general solution of equation (\ref{fe2}) as follows
\[
y(t)=y(a)+(y(b)-y(a))(2\alpha-1)\cdot (\frac{b-t}{b-a})^{2\alpha-1}\int_{1}^{\frac{b-a}{b-t}} (w-1)^{\alpha-1}\, w^{\alpha -1}\, dw,
\]
in which $y(a)$ and $y(b)$ are constants to be determined by imposing the boundary conditions. 
\begin{fe}
For $ 0<\alpha<1$ consider the following two term fractional Sturm-Liouville equation
\begin{equation}\label{fe3}
-^{c}\mathcal{D}_{0^+}^{\alpha}\circ   \mathcal{D}^{\alpha}_{0^+}y(t)=\lambda y(t),\qquad t > 0.
\end{equation}
\end{fe}
 By taking the Laplace transform on both sides of  equation (\ref{fe3}) and using Property~\ref{lc}, we have
\[
-s^{\alpha}\mathfrak{L}\{\mathcal{D}^{\alpha}_{0^+}y(t)\}+s^{\alpha-1}\mathcal{D}^{\alpha}_{0^+}y(t)|_{t=0}=\lambda\mathfrak{L}\{y(t)\},
\]
and then using Property~\ref{lrl}, we get
\[
-s^{\alpha}\left(s^{\alpha}\mathfrak{L}\{y(t)\}-\mathcal{I}_{0^+}^{1-\alpha}y(t)|_{t=0}\right)+s^{\alpha-1}\mathcal{D}_{0^+}^{\alpha}y(t)|_{t=0}=\lambda\mathfrak{L}\{y(t)\}.
\]
It can be readily obtained 
\[
\mathfrak{L}\{y(t)\}=\frac{s^{\alpha}}{\lambda+s^{2\alpha}}\mathcal{I}^{1-\alpha}_{0^+}y(t)|_{t=0}+\frac{s^{\alpha-1}}{\lambda+s^{2\alpha}}\mathcal{D}^{\alpha}_{0^+}y(t)|_{t=0}.
\]
By taking the inverse Laplace transform in order to gain $y(t)$, we get
\begin{equation}\label{lg}
y(t)=\mathfrak{L}^{-1}\{\frac{s^{\alpha}}{\lambda+s^{2\alpha}}\}\cdot\mathcal{I}^{1-\alpha}_{0^+}y(t)|_{t=0}+\mathfrak{L}^{-1}\{\frac{s^{\alpha-1}}{\lambda+s^{2\alpha}}\}\cdot\mathcal{D}^{\alpha}_{0^+}y(t)|_{t=0}.
\end{equation}
On the other hand,
\begin{equation}{\nonumber}
\begin{split}
\mathfrak{L}^{-1}\{\frac{s^{\alpha}}{\lambda+s^{2\alpha}}\}&=\mathfrak{L}^{-1}\{\frac{s^{\alpha}}{\lambda}\frac{1}{1+\left(\frac{s}{\lambda^{\frac{1}{2\alpha}}}\right)^{2\alpha}}\}\\
&=\mathfrak{L}^{-1}\{\frac{s^{\alpha}}{\lambda}-\frac{s^{3\alpha}}{\lambda^2}+\frac{s^{5\alpha}}{\lambda^3}-\cdots\}\qquad \text{as}\quad|\lambda|\rightarrow\infty\\
&=\mathfrak{L}^{-1}\{\sum_{k=1}^{\infty}(-1)^{k-1}\frac{s^{(2k-1)\alpha}}{\lambda^k}\}\\
&=\sum_{k=1}^{\infty}\frac{(-1)^{k-1}}{\lambda^kt^{(2k-1)\alpha+1}\Gamma(-(2k-1)\alpha)}\qquad \text{(by virtue of Property~\ref{ils})}\\
&=-\frac{1}{t^{1-\alpha}}\sum_{k=1}^{\infty}\left(\frac{-1}{\lambda t^{2\alpha}}\right)^k\frac{1}{\Gamma(\alpha-2\alpha k)}.
\end{split}
\end{equation}
Consequently, assuming $z=-\lambda t^{2\alpha}$, $\delta=2\alpha$ and $\theta=\alpha$ in (\ref{asy}), we can write from above equation
\[
\mathfrak{L}^{-1}\{\frac{s^{\alpha}}{\lambda+s^{2\alpha}}\}=\frac{1}{t^{1-\alpha}}\left(E_{2\alpha,\alpha}(-\lambda t^{2\alpha})-\frac{1}{2\alpha}(-\lambda t^{2\alpha})^{\frac{1-\alpha}{2\alpha}}\exp((-\lambda t^{2\alpha})^{\frac{1}{2\alpha}})\right),
\]
as $|\lambda|\rightarrow\infty$. For real positive eigenvalues $\lambda$ and in accordance with $t\in [0,1]$, we have $\arg(-\lambda t^{2\alpha})=\pi$. Hence, the second term of the right hand side of the above equation vanishes due to Property~\ref{asym} and we get
\begin{equation}\label{mi1}
\mathfrak{L}^{-1}\{\frac{s^{\alpha}}{\lambda+s^{2\alpha}}\}=\frac{1}{t^{1-\alpha}}E_{2\alpha,\alpha}(-\lambda t^{2\alpha}).
\end{equation}
In similar manner, we can obtain 
\begin{equation}\label{mi2}
\mathfrak{L}^{-1}\{\frac{s^{\alpha-1}}{\lambda+s^{2\alpha}}\}=t^{\alpha}E_{2\alpha,\alpha+1}(-\lambda t^{2\alpha}).
\end{equation}
So, using (\ref{lg}), (\ref{mi1}) and (\ref{mi2}), we can write the general solution of (\ref{fe3})  as,
\begin{equation}\label{sf3}
y(t)=\mathcal{I}_{0^+}^{1-\alpha}y(t)|_{t=0}\cdot \frac{1}{t^{1-\alpha}}E_{2\alpha,\alpha}(-\lambda t^{2\alpha})+ \mathcal{D}_{0^+}^{\alpha}y(t)|_{t=0}\cdot t^{\alpha}E_{2\alpha,\alpha+1}(-\lambda t^{2\alpha}),
\end{equation}
in which $\mathcal{I}_{0^+}^{1-\alpha}y(t)|_{t=0}$ and $\mathcal{D}_{0^+}^{\alpha}y(t)|_{t=0}$ are constants to be determined using the initial conditions. The two linearly independent functions $y_1(t;\alpha,\lambda)=\frac{1}{t^{1-\alpha}}E_{2\alpha,\alpha}(-\lambda t^{2\alpha})$ and $y_2(t;\alpha,\lambda)=t^{\alpha}E_{2\alpha,\alpha+1}(-\lambda t^{2\alpha})$ satisfy the equation (\ref{fe3}), separately so, they form a fundamental set of solutions. It is readily verified that $y_2^\prime (t;\alpha,\lambda))=y_1(t;\alpha,\lambda)
$. Note that \eqref{sf3} may have been obtained alternately using the method of successive approximations.
\begin{remark}
When $\alpha$ approaches $1$, equation (\ref{fe3}) turns into $-y^{\prime\prime}=\lambda y$ and its fundamental set of solutions, i.e., $\{t^{\alpha-1}E_{2\alpha,\alpha}(-\lambda t^{2\alpha}),t^{\alpha}E_{2\alpha,\alpha+1}(-\lambda t^{2\alpha})\}$ reduces to $\{\cos(\sqrt{\lambda}t),\frac{\sin(\sqrt{\lambda}t)}{\sqrt{\lambda}}\}$ due to (\ref{e1}), (\ref{e2}) and (\ref{pmi}). This  shows that our results are a  generalization of classical the ones.
\end{remark}

\section{Eigenvalues of our Fractional Sturm-Liouville problem}
Consider (\ref{fe3}) on $[0,1]$ with boundary conditions
\begin{equation}\label{fcon}
\mathcal{I}_{0^+}^{1-\alpha}y(t)|_{t=0}=0,\quad \text{and}\quad \mathcal{I}_{0^+}^{1-\alpha}y(t)|_{t=1}=0.
\end{equation}
Imposing these boundary conditions on the general solution of (\ref{sf3}) with $\mathcal{D}_{0^+}^{\alpha}y(t)|_{t=0} \neq 0$  and using that
\[
\mathcal{I}_{0^+}^{1-\alpha}\{t^{\alpha}E_{2\alpha,\alpha+1}(-\lambda t^{2\alpha})\}=tE_{2\alpha,2}(-\lambda t^{2\alpha}),
\]
we get
\begin{equation}\label{eigeq}
E_{2\alpha,2}(-\lambda)=0,
\end{equation}
which defines our characteristic equation for the eigenvalues. We note that $E_{2\alpha,2}(-\lambda)$ is an entire function of $\lambda$ order $1/2\alpha$ and type 1, [\cite{kil}, p.42], and since it must be of fractional order for all $ 1/2 < \alpha < 1$ by classical complex analysis it must have infinitely many zeros for all such $\alpha$.

Now the Mittag-Leffler function $E_{2\alpha,2}(-\lambda)$ for $1/2<\alpha<1$ can be decomposed into two parts \cite{gor} given by
\begin{equation}\label{plus}
\lambda^{\frac{1}{2\alpha}}E_{2\alpha,2}(-\lambda)=f_{2\alpha,2}(-\lambda)+g_{2\alpha,2}(-\lambda)
\end{equation}
where
\begin{equation}\label{f}
f_{2\alpha,2}(-\lambda)=\int_0^{\infty}e^{-r\lambda^{\frac{1}{2\alpha}}}k_{2\alpha,2}(r)dr
\end{equation}
with
\[
k_{2\alpha,2}(r)=\frac{1}{\pi}\frac{r^{2\alpha-2}(-\sin(2\alpha \pi))}{r^{4\alpha}+2r^{2\alpha}\cos(2\alpha\pi)+1}
\]
and
\begin{equation}\label{g}
g_{2\alpha,2}(-\lambda)=\frac{1}{\alpha}e^{\lambda^{\frac{1}{2\alpha}}\cdot\cos(\frac{\pi}{2\alpha})}\cdot \cos\left(\lambda^{\frac{1}{2\alpha}}\sin(\frac{\pi}{2\alpha})-\frac{\pi}{2\alpha}\right).
\end{equation}
It is obvious that 
\[
g_k=\left(\frac{(k+\frac{1}{2}+\frac{1}{2\alpha})\pi}{\sin(\frac{\pi}{2\alpha})}\right)^{2\alpha},\qquad k=-1,0,1,\cdots,
\]
are positive zeros of $g_{2\alpha,2}(-\lambda)$ and $g_{2\alpha,2}(0)=\frac{1}{\alpha}\cos\frac{\pi}{2\alpha}$. We note that this function exhibits oscillations with an amplitude which decays  exponentially, since $\cos(\frac{\pi}{2\alpha})$ is negative as long as $1/2<\alpha<1$.

The function $k_{2\alpha,2}(r)$ is always positive due to the fact that $\sin(2\alpha\pi)$ is negative for $1/2<\alpha<1$ and the denominator is greater than $(r^{2\alpha}-1)^2$ which is nonnegative for such $\alpha$. Also, we have
\[
f_{2\alpha,2}(0)=\int_0^{\infty}k_{2\alpha,2}(r)dr=-\frac{1}{\alpha}\cos(\frac{\pi}{2\alpha})>0
\]
for $1/2<\alpha<1$ and according to Riemann-Lebesgue Lemma \cite{leb}, the function $f_{2\alpha,2}(-\lambda)$ approaches zero  asymptotically as $\lambda\rightarrow\infty$. Recall that the form of $E_{2\alpha, 2}(-\lambda)$ (cf., \eqref{mittag2}) shows that as $\lambda \to \infty$ its zeros as a function of $\lambda$, if any, must lie along the negative axis. Its asymptotic form as a function of a complex variable can be found in [\cite{erdelyi}, p. 210, eq. (22)] with the obvious substitutions. Specifically, we know that
$$E_{2\alpha, 2}(z) =\frac{1}{2\alpha} \sum_{m} t_m^{-1}\, e^{t_m} - O\left (|\lambda|^{-(N-1)}\right ),\quad {\rm as}\quad \lambda \to \infty.$$
Here $t_m = z^{1/2\alpha}\, e^{\pi i m/\alpha}$, $m$ is an integer,  and $-2\alpha \pi < \arg z +2\pi m < 2\alpha \pi$. One can show directly from this that, for each $\lambda > 0$,  $$E_{2\alpha, 2}(-\lambda) \to \frac{\sin\sqrt{\lambda}}{\sqrt{\lambda}}$$ as $\alpha \to 1$, as expected in the classical case. The latter display on the right gives the dispersion relation for the eigenvalues (all positive) of the Dirichlet problem for $-y^{\prime\prime}=\lambda\, y$ on $[0,1]$. It also follows that the eigenvalues of our problem approach those of the Dirichlet problem for $-y^{\prime\prime}=\lambda\, y$, as the zeros of $E_{2\alpha, 2}(-\lambda)$ approach those of $\frac{\sin\sqrt{\lambda}}{\sqrt{\lambda}}$

Since $E_{2\alpha, 2}(-\lambda)$ is an entire function of (fractional) finite order $1/2\alpha$ [\cite{erdelyi}, p. 208], we can ascertain that, for each $\alpha$, there must be a {\it finite} number of zeros along the negative $\lambda$ axis whenever $1/2 <\alpha < 1$. In fact, for each such $\alpha$, $E_{2\alpha, 2}(-\lambda) > 0$ for all sufficiently large $\lambda$, where the $\lambda$ generally depends on $\alpha$. This now implies that 
\begin{equation}\label{xine}
f_{2\alpha,2}(-\lambda) > |g_{2\alpha,2}(-\lambda)| 
\end{equation}
for all sufficiently large (positive) $\lambda$. Next, in the remaining {\it finite}  $\lambda$-interval, we may have subintervals wherein 
\begin{equation}\label{ine}
|g_{2\alpha,2}(-\lambda)|\geq f_{2\alpha,2}(-\lambda).
\end{equation}
In such an exceptional interval in which $g_{2\alpha,2}(-\lambda)<0$, call it $I$, \eqref{g} implies that $ f_{2\alpha,2}(-\lambda) = - g_{2\alpha,2}(-\lambda)$ at least twice. Hence $E_{2\alpha, 2}(-\lambda)$ has two zeros in $I$ (each of which is an eigenvalue).

Denote the totality of such intervals in which $g_{2\alpha,2}(-\lambda)<0$ by $I_n$, $n=0, 1, 2, \ldots, N^*-1$, where $N^*$ depends on $\alpha$. Its value will be specified below. A glance at (\ref{g}) shows that the $n^{th}$ such interval $I_n$ is given by the set of all $\lambda$ such that,
\[
(4n+1)\frac{\pi}{2}<\lambda^{\frac{1}{2\alpha}}\sin(\frac{\pi}{2\alpha})-\frac{\pi}{2\alpha}<(4n+3)\frac{\pi}{2},\quad n=0,1,\cdots.
\]
In other words, 
\begin{equation}\label{neg}
I_n:=\left(\left(\frac{(2n+\frac{1}{2}+\frac{1}{2\alpha})\pi}{\sin(\frac{\pi}{2\alpha})}\right)^{2\alpha},\left(\frac{(2n+\frac{3}{2}+\frac{1}{2\alpha})\pi}{\sin(\frac{\pi}{2\alpha})}\right)^{2\alpha}\right),\quad n=0,1,\cdots.
\end{equation}
Recall, however, that \eqref{xine} must persist for all large $\lambda$ so that, insofar as eigenvalues are concerned, this list must be necessarily {\it finite} for each $\alpha$ from previous considerations and consists of the intervals $I_n$, $n=0, 1, 2, \ldots, N^*-1$. We conclude that each interval $I_n$ under consideration contains two eigenvalues and their total number is then $2N^*$. 
\subsection{Estimating the number of real eigenvalues}
It is worth noting again that the number of real eigenvalues is finite while the number of non-real eigenvalues is infinite (as the entire function is of fractional order for our range of $\alpha$). We do not, however, address the non-real spectrum of our problem and, in the following, consider the real spectrum exclusively.

Determining the extreme points of $g_{2\alpha,2}(-\lambda)$ we find,
\[
g'_{2\alpha,2}(-\lambda)=\frac{e^{\lambda^{\frac{1}{2\alpha}}\cos(\frac{\pi}{2\alpha})}\cdot \lambda^{\frac{1}{2\alpha}-1}}{2\alpha^2}\left(\cos(\frac{\pi}{2\alpha})\cos(\lambda^{\frac{1}{2\alpha}}\sin(\frac{\pi}{2\alpha})-\frac{\pi}{2\alpha})-\sin(\frac{\pi}{2\alpha})\sin(\lambda^{\frac{1}{2\alpha}}\sin(\frac{\pi}{2\alpha})-\frac{\pi}{2\alpha})\right)=0
\]
which yields
\[
g'_{2\alpha,2}(-\lambda)=\frac{e^{\lambda^{\frac{1}{2\alpha}}\cos(\frac{\pi}{2\alpha})}\cdot \lambda^{\frac{1}{2\alpha}-1}}{2\alpha^2}\left(\cos(\lambda^{\frac{1}{2\alpha}}\sin(\frac{\pi}{2\alpha})\right)=0.
\]
Therefore, 
\[
z_k=\left(\frac{(k+\frac{1}{2})\pi}{\sin(\frac{\pi}{2\alpha})}\right)^{2\alpha}\qquad, k=0,1,\cdots,
\]
are the extreme points of $g_{2\alpha,2}(-\lambda)$ and at these points we have
\begin{equation}\nonumber
g_{2\alpha,2}(-z_k)=\frac{1}{\alpha}e^{(k+\frac{1}{2})\pi\cdot \cot(\frac{\pi}{2\alpha})}\cdot(-1)^k\sin(\frac{\pi}{2\alpha})\quad \left\{
\begin{array}{ll}
>0&
k \quad\text{even}, \\ \\
<0 & k \quad\text{odd},
\end{array}\right.
\end{equation}
for $1/2<\alpha<1$. Since we are interested in that part of $g_{2\alpha,2}(-\lambda)$ where it is negative, we need to look at those extreme points arising from an odd value of $k$, that is,
\[
z^{odd}_k=\left(\frac{(2k+\frac{3}{2})\pi}{\sin(\frac{\pi}{2\alpha})}\right)^{2\alpha},\qquad k=0,1,\cdots.
\]
The above considerations imply that for each fixed $1/2<\alpha<1$, there exists an $N^{\ast}(\alpha)\in \mathbb{N}$ such that for all $N\geq N^{\ast}$ we have
\begin{equation}\label{ine2}
\left|g_{2\alpha,2}(-z_N^{odd})\right|<f_{2\alpha,2}(-z_N^{odd}),
\end{equation}
and such that there is no zero of $E_{2\alpha,2}(-\lambda)$ in these negative intervals of $g_{2\alpha,2}(-\lambda)$, namely $I_n$ for $n=N^{\ast},N^{\ast}+1,\ldots$. The inequality (\ref{ine2}) can be rewritten as follows,
\begin{equation}\label{ine3}
\frac{1}{\alpha}e^{(2N^{\ast}+\frac{3}{2})\pi\cdot \cot(\frac{\pi}{2\alpha})}\cdot\sin(\frac{\pi}{2\alpha})<\frac{1}{\pi}\int_0^{\infty}e^{-r\frac{(2N^{\ast}+\frac{3}{2})\pi}{\sin(\frac{\pi}{2\alpha})}}\frac{r^{2\alpha-2}(-\sin(2\alpha \pi))}{r^{4\alpha}+2r^{2\alpha}\cos(2\alpha\pi)+1}dr.
\end{equation}
From this, for given $\alpha$, we can solve for for $N^{\ast}$ implicitly and so numerically. Since $f_{2\alpha,2}$ is decreasing and $g_{2\alpha,2}$ is an oscillating cosine function with an exponentially decaying amplitude (as $1/2 < \alpha < 1$), there are two zeros of $E_{2\alpha,2}(-\lambda)$ in each interval $I_n$ for $n=0,1,\cdots,N^{\ast}-1$ in which $g_{2\alpha,2} < 0$. 

{\bf Remark:}\quad It must be mentioned that the numerical solution of the integral appearing in the right hand side of (\ref{ine3})  for fixed $\alpha$ can be obtained by using various library packages such as \emph{NIntegrate} in \emph{Mathematica} or even \emph{Maple}.

Since the intervals where $g_{2\alpha,2} < 0$ are given by \eqref{neg}, we can readily deduce the following two results.
\begin{prop}
The $2\alpha$-th root of the two first eigenvalues of  (\ref{fe3}) - (\ref{fcon}), namely $\rho_{1}$ and $\rho_{2}$ if they exist, are in the interval  
\begin{equation}\label{sint}
\tilde{I}_{0}=\left(\left(\frac{\frac{1}{2}+\frac{1}{2\alpha})\pi}{\sin(\frac{\pi}{2\alpha})}\right),\left(\frac{(\frac{3}{2}+\frac{1}{2\alpha})\pi}{\sin(\frac{\pi}{2\alpha})}\right)\right),
\end{equation}
and $\lim_{\alpha\rightarrow 1} |\tilde{I}_{0}|=\pi$, which corresponds to the classical case.
\end{prop}
\begin{prop}
The $2\alpha$-th root of the two largest eigenvalues denoted by  $\rho_{2N^{\ast}-1}$ and $\rho_{2N^{\ast}}$, are in the interval 
\begin{equation}\label{bint}
\tilde{I}_{N^{\ast}-1}=\left(\left(\frac{(2N^{\ast}-\frac{3}{2}+\frac{1}{2\alpha})\pi}{\sin(\frac{\pi}{2\alpha})}\right),\left(\frac{(2N^{\ast}-\frac{1}{2}+\frac{1}{2\alpha})\pi}{\sin(\frac{\pi}{2\alpha})}\right)\right),
\end{equation}
and $\lim_{\alpha\rightarrow 1} |\tilde{I}_{N^{\ast}-1}|=\pi$.
\end{prop}
Table~\ref{tab1} gives numerical results for different values of $1/2<\alpha<1$.  Using  (\ref{ine3}), (\ref{sint}), and  (\ref{bint})  we can then find the number of eigenvalues of (\ref{fe3}) -  (\ref{fcon}). The intervals containing the $2\alpha$-th root of the smallest and largest eigenvalue are shown in the two columns on the right.
 \begin{table}[H] 
\caption{Numerical results for many different values of $\alpha$. The intervals containing the $2\alpha$-th root of the smallest and largest eigenvalue are shown in the two columns on the right.} 
 \centering 
 \begin{tabular}{c c c c  } 
 \hline\hline  \\ [-2ex]
 $\alpha$ (Fractional order) & $2N^{\ast}$ (Number of eigenvalues) & $\tilde{I}_0$ &$\tilde{I}_{N^{\ast}-1}$ \\ [0.5ex] 
 \hline 
 0.78 & 0              & -----& ----- \\  
 0.80 & 2         &(3.82549,7.22593)&(3.82549,7.22593)\\ 
 0.82&  2       &(3.70445,7.04252)&(3.70445,7.04252)\\ 
 0.84&2   &  (3.60076,6.88842) &(3.60076,6.88842) \\
 0.86&4 & (3.51148,6.75866) & (10.0058,13.253) \\
 0.88& 4 &(3.43428,6.64934) & (9.86441,13.0795)  \\
 0.90& 8 & (3.36728,6.55734) & (22.5076,25.6977) \\
 0.92& 10 & (3.309,6.48013) & (28.678,31.8492)  \\
  0.94& 18 & (3.25822,6.41567) & (53.7774,56.9349)  \\
  0.96&32   &  (3.21392,6.36226) &(97.6639,100.812) \\
 0.98&84 & (3.17528,6.31849) & (260.918,264.062) \\
 0.981& 90 &(3.17348,6.31653) & (279.762,282.905)  \\
 0.982& 98 & (3.1717,6.3146) & (304.89,308.033) \\
 0.983& 104 & (3.16993,6.31268) & (323.731,326.873)  \\
  0.984& 114 & (3.16817,6.31079) & (355.141,358.284)  \\
   0.985& 124 & (3.16642,6.30891) & (386.55,389.693) \\
 0.989& 182 & (3.15955,6.30162) & (568.733,571.875)  \\
  0.9898& 200 & (3.15819,6.3002) & (625.275,628.417)  \\  
 [1ex] 
 \hline 
 \end{tabular} 
 \label{tab1} 
 \end{table}

{\bf Remark:}\quad Observe that the length of these intervals, $\tilde{I}_0$ and $\tilde{I}_{N^{\ast}-1}$, approaches $\pi$ as $\alpha \to 1^-$ (as expected in the classical case).  In addition, a glance at \eqref{neg} shows that the asymptotic behavior of an  eigenvalue, let's call it $\lambda_n (\alpha)$, contained in \eqref{neg}, for fixed $n$, is given by
$$\lambda_n (\alpha) \sim \left (\frac{n\pi}{\sin (\frac{\pi}{2\alpha})}\right)^{2\alpha},$$
as $\alpha \to 1^-$. This corresponds exactly with the well known classical asymptotic estimate $\lambda_n \sim n^2\,\pi^2$ as $n \to \infty$.

 \section{Conclusion} 
In this article, we gave an introduction to three distinct fractional Sturm-Liouville equations consisting of  various compositions of a Riemann-Liouville and Caputo fractional derivatives simultaneously operating on an unknown function. It must be mentioned that the potential function is absent in these equations for which the fundamental set of solutions have been derived.

Of these operators we distinguished one of them and solved a Dirichlet type problem completely all the while showing that, in the limit as $\alpha\to 1^-$ we recover the two-term Sturm-Liouville problem in the case of Dirichlet boundary conditions.
 
Our results might set the groundwork for many semi-analytical and numerical methods such as Homotopy Perturbation Method (HPM) and Variational Iteration Method (VIM)  once a potential function is introduced, through which the eigenvalues and the eigenfunctions of a Sturm-Liouville  problem are found via the fundamental  solutions of the linear part of the equation. It is also possible to implement the classical {\it Variation of Parameters} method to obtain the general solution of a fractional Sturm-Liouville equation analytically in the presence of a  potential function though we have not done this here.



\end{document}